\documentclass[journal,12pt,onecolumn,draftclsnofoot,]{IEEEtran}

\pdfoutput=1
\IEEEoverridecommandlockouts                              % This command is only needed if 
\newtheorem{Remark}{Remark}
\newtheorem{Corollary}{Corollary}

\newenvironment{Proof}{\noindent{\em Proof:\/}}{\hfill $\Box$\par}
\newtheorem{Theorem}{Theorem}
\newtheorem{Lemma}{Lemma}

\newtheorem{Assumption}{Assumption}

\usepackage{graphicx}
\usepackage{amsmath,bm}
\usepackage{amssymb}
\usepackage{algorithm}
\usepackage{algpseudocode}
\usepackage{cite}
\usepackage{hyperref}
\usepackage[caption=false]{subfig}
\hypersetup{hidelinks}
\algdef{SE}[DOWHILE]{Do}{doWhile}{\algorithmicdo}[1]{\algorithmicwhile\ #1}%

\newcommand{\mathactivatecomma}{%
  \begingroup\lccode`~=`\,
  \lowercase{\endgroup\edef~}{\mathchar\the\mathcode`\,\penalty0 }}

\algnewcommand{\Initialize}[1]{%
  \State \textbf{Initialize: $i \in \mathcal{V}$}
  \Statex \hspace*{\algorithmicindent}\parbox[t]{.8\linewidth}{\raggedright #1}
}

\algnewcommand{\Iteration}[1]{%
  \State \textbf{Iteration $(t\geq 0)$: $i \in \mathcal{V}$}
  \Statex \hspace*{\algorithmicindent}\parbox[t]{.8\linewidth}{\raggedright #1}
}

\algnewcommand{\Output}[1]{%
  \State \textbf{Output: $i \in \mathcal{V}$}
  \Statex \hspace*{\algorithmicindent}\parbox[t]{.8\linewidth}{\raggedright #1}
}

\begin{document}

\bstctlcite{IEEEexample:BSTcontrol}

\title{Randomized Gradient-Free Distributed Online Optimization via a Dynamic Regret Analysis}

\author{Yipeng Pang and Guoqiang Hu % <-this % stops a space
% \thanks{This research was supported by Singapore Ministry of Education Academic Research Fund Tier 1 RG180/17(2017-T1-002-158).}% <-this % stops a space
\thanks{Y. Pang and G. Hu are with the School of Electrical and Electronic Engineering, Nanyang
Technological University, 639798, Singapore
        { ypang005@e.ntu.edu.sg, gqhu@ntu.edu.sg}.}%
}

\maketitle
% \thispagestyle{empty}
% \pagestyle{empty}

% \begin{keyword}                           % Five to ten keywords,  
% Distributed optimization; multi-agent system; gradient-free optimization.               % chosen from the IFAC 
% \end{keyword}                             % keyword list or with the 
                                          % help of the Automatica 
                                          % keyword wizard

\begin{abstract}                          % Abstract of not more than 200 words.
This work considered an online distributed optimization problem, with a group of agents whose local objective functions vary with time. Moreover, the value of the objective function is revealed to the corresponding agent after the decision is executed per time-step. Thus, each agent can only update the decision variable based on the revealed value and information collected from the neighbors, without the knowledge on the explicit expression of the objective function. To solve this problem,
% Followed by our preliminary work in \cite{Pang2019b}, 
an online gradient-free distributed projected gradient descent (DPGD) algorithm is presented, where each agent locally approximates the gradient based on two point values. 
% To analyze the convergence property, 
% the notion of the dynamic regret is adopted to measure the difference between the total costs incurred by the agent's state trajectory and the offline centralized optimal solution where the objective functions are available {\it a priori}.
With some standard assumptions on the communication graph and the objective functions, we provide the bound for the dynamic regret as a function of the minimizer path length, step-size and smoothing parameter. Under appropriate selections of the step-size and smoothing parameter, we prove that the dynamic regret is sublinear with respect to the time duration $T$ if the minimizer path length also grows sublinearly.
Finally, the effectiveness of the proposed algorithm is illustrated through numerical simulations.
\end{abstract}

% \begin{abstract}
% In this note, we study an online multi-agent optimization problem where the objective functions of agents vary with time. A gradient-free distributed algorithm is proposed
% \end{abstract}


\section{Introduction}
Many applications related to the coordination among multiple agents can be formulated as a distributed optimization problem, where the global objective function is often cast as a summation of each individual's local objective function. In this problem, the agent updates its own decision variable based on the local information and the information from its neighbors in the network, to achieve the convergence to the minimizer of the global objective function. 
Such problem commonly exists in a variety of applications, including parameter estimation, source localization, utility maximization, resource allocation, path-planning, \textit{etc}. 
In some situations, the uncertainties in the environment may influence the objective functions, which in turn, affect the performance of the designed algorithms. One approach to cater such issues is through stochastic methods, which has been studied in \cite{Mertikopoulos2017, Agarwal2012, Ram2010a}. However, in many practical cases, especially when mobile agents are involved, the optimization problem is often in a highly dynamic environment, resulting in a time-varying objective function. Thus, the techniques used in fixed or distributed static optimization cannot be formally applied, leading to the study of online optimization framework, where the objective functions assigned to agents vary with time and these variations are revealed to agents only in hindsight.

\textit{Related Works:} Vast results on distributed online optimization have been reported in \cite{Sun2017,Yan2013,Akbari2017,Hosseini2013,Hosseini2016,Mateos-Nunez2014,Lee2017}, where the considered regret is
% The work in \cite{Sun2017} proposed a continuous-time gradient-based algorithm for time-varying quadratic problems. For the problems with strongly convex objective functions, a distributed online learning based algorithm was introduced in \cite{Yan2013}, which shows a regret bound of $\mathcal{O}(\ln{T})$. The work in \cite{Akbari2017} adopted a push-sum based algorithm to remove the doubly-stochastic condition on the weighting matrix, obtaining a regret bound of $\mathcal{O}((\ln{T})^2)$. For problems with general convex objective functions, the works in \cite{Hosseini2013} and \cite{Hosseini2016} proposed distributed dual-subgradient averaging algorithms, deducing a regret bound of $\mathcal{O}(\sqrt{T})$. This result was recovered in \cite{Mateos-Nunez2014} by subgradient descent algorithms with proportional-integral disagreement feedback, and in \cite{Lee2017} by Nesterov's dual averaging methods.
% It is worth noting that the regret used in the aforementioned literature is 
{\it static regret}, defined by the difference between the total incurred cost and the cost of the best fixed decision in hindsight. On the other hand, to study the scenario where the functions and the decision variables evolve simultaneously instead of a single best fixed decision, the notion of {\it dynamic regret} was brought forward to characterize how much one regrets working in an online setting as opposed to the offline solution with full knowledge of past and future observations. 
However, achieving a sublinear dynamic regret is usually not possible in the worst case scenario, and the regret is usually characterized with some problem regularities, such as the minimizer path length $\omega_T$.
With this concept, the works in \cite{Shahrampour2017,Shahrampour2018} proposed decentralized mirror descent methods for the problem with unknown unstructured noise, and obtained a regret bound of $\mathcal{O}(\sqrt{\omega_TT})$. In \cite{Yi2020}, time-varying coupled inequality constraints were considered and a primal-dual dynamic mirror descent algorithm was designed to achieve a regret of $O(\max\{\omega_TT^\kappa,T^{\max\{\kappa,1-\kappa\}}\})$. The work in \cite{Lu2020} developed a distributed algorithm based on an auxiliary optimization strategy for pseudo-convex objective functions and showed a regret bound of $\mathcal{O}(\sqrt{\omega_T}T^{3/4}\ln{T})$. The work in \cite{Zhang2020} employed gradient-tracking techniques for problems with strongly convex objective functions and achieved a regret bound of $O(C+\omega_T+\nu_T)$ with a constant step-size, where $C$ is related to the initial conditions and $\nu_T$ is the gradient path length. 
It should be noted that all these methods make use of the derivatives in the process of the optimization, which implicitly presumes that the derivatives can be obtained directly. However, there are many applications where the gradient information is not available to use, then these methods cannot be applied. This inspires the study of gradient-free techinques, which has been studied in our recent works \cite{Pang2017,Pang2020,Pang2018} for distributed static optimization problems. 
% {\it e.g.}, \cite{Nesterov2017,Yuan2015,Li2015,Chen2017,Yuan2015a,Pang2017,Pang2020,Pang2018}, to list a few. A Gaussian-smoothing technique was firstly proposed in \cite{Nesterov2017} to solve the general unconstrained optimization problem. This technique was applied to a distributed optimization problem in \cite{Yuan2015,Li2015} and further improved in \cite{Chen2017,Yuan2015a}.  
% The variants of this technique were also adopted in \cite{Pang2017,Pang2020,Pang2018} to solve the distributed static optimization problems. 
For distributed online optimization problems, the works in \cite{Lee2017,Shahrampour2017,Shahrampour2018} considered the scenarios where the gradient of the objective function is coupled with noise, hence proposed stochastic gradient methods. However, gradient-free optimization schemes have received little attention in online optimization problems in general. The only relevant work is \cite{Shames2020}, which utilized finite differences to estimate the derivative, but is not distributed.

\textit{Contributions:} In this work, we propose an online gradient-free optimization algorithm based on Gaussian smoothing technique \cite{Nesterov2017} to locally approximate the function derivative. Some preliminary results of this work have been presented in \cite{Pang2019b}, where static regret was adopted in the analysis. Compared with the existing literature including our work in \cite{Pang2019b}, the major contributions of this paper are threefold. 1). Different from the consensus-based approaches in most existing distributed online optimization literature \cite{Yan2013,Hosseini2013,Hosseini2016,Mateos-Nunez2014,Shahrampour2017,Shahrampour2018,Yi2020,Lu2020,Zhang2020}, the proposed algorithm adopted the surplus-based method, which was firstly proposed in \cite{Cai2012} to solve average consensus problems, to remove the doubly-stochastic requirement on the weighting matrix\footnote{We note that there are other consensus techniques to remove the doubly-stochastic requirement on the weighting matrix, \textit{e.g.}, push-sum methods, which have been studied in distributed static optimization, see \cite{Benezit2010,Nedic2015,Yuan2015a}.}, and hence enables the implementation in any strongly connected and fixed digraphs.
2). Different from \cite{Yan2013,Akbari2017,Hosseini2013,Hosseini2016,Mateos-Nunez2014,Lee2017} where static regret was utilized to analyze the performance of the algorithm, this paper adopts the notion of the dynamic regret to study the scenario where the functions and the decision variables evolve simultaneously instead of a single best fixed decision.
3). In contrast to the online optimization literature with dynamic regret \cite{Shahrampour2017,Shahrampour2018,Yi2020,Lu2020,Zhang2020} where the gradient information can be obtained directly, we considered gradient-free settings and adopted the Gaussian smoothing technique \cite{Nesterov2017} to locally approximate the function derivative. The same technique was adopted in our work \cite{Pang2019b} to achieve a sublinear static regret plus a linear error term due to the constant smoothing parameter. In this work, we eliminate the linear error term by choosing a non-increasing smoothing parameter sequence. Specifically, we provide the bound for the dynamic regret as a function of the minimizer path length, step-size and smoothing parameter. Under appropriate selections of the step-size and smoothing parameter, we prove that the dynamic regret is sublinear if the minimizer path length also grows sublinearly.

% In the subsequent sections, the notations used in this paper is introduced in section~\ref{sec:problem_formulation}, followed by the formal definition of the problem. The proposed algorithm is detailed in section~\ref{sec:distr_opt}, and the performance analysis via dynamic regret is conducted in section~\ref{sec:analysis}. We verify the effectiveness of the proposed algorithm by some numerical simulations in section~\ref{sec:simulation}. Section~\ref{sec:conclusion} concludes the paper.

\textit{Notations:} We use $\mathbb{R}^p$ to denote the set of $p$-dimensional column vectors, and $\mathbf{E}[a]$ for the expected value of a random variable $a$. For any two vectors $x$ and $y$, the operator $\langle x,y\rangle$ denotes the inner product of $x$ and $y$, and $\Pi_{\Omega}[x]$ denotes the projection of $x$ on a set $\Omega$, \textit{i.e.}, $\Pi_{\Omega}[x] = \arg\min_{\hat{x} \in \Omega} \|\hat{x} - x\|^2$. $\|x\|$ is the standard Euclidean norm of $x$, \textit{i.e.}, $\|x\| = \sqrt{x^\top x}$. For a set $\mathcal{N}$, $|\mathcal{N}|$ represents the number of elements in $\mathcal{N}$.

%%%%%%%%%%%%%%%%%%%%%%%%%%%%%%%%%%%%%%%%%%%%%%%%%%%%%%%%%%%%%%%%%%%%%%%%%%%%%%%%
\section{Problem Formulation and Preliminaries}\label{sec:problem_formulation}

% This section introduces the notations used. Then, the problem is formally defined, followed by some preliminary results.

% \subsection{Notations}
% We use $\mathbb{R}^p$ to denote the set of $p$-dimensional column vectors, and $\mathbf{E}[a]$ for the expected value of a random variable $a$. For any two vectors $x$ and $y$, the operator $\langle x,y\rangle$ denotes the inner product of $x$ and $y$, and $\Pi_{\Omega}[x]$ denotes the projection of $x$ on a set $\Omega$, \textit{i.e.,} $\Pi_{\Omega}[x] = \arg\min_{\hat{x} \in \Omega} \|\hat{x} - x\|^2$. $\|x\|$ is the standard Euclidean norm of $x$, \textit{i.e.}, $\|x\| = \sqrt{x^\top x}$. For a set $\mathcal{N}$, $|\mathcal{N}|$ represents the number of elements in $\mathcal{N}$.

\subsection{Problem Definition}
Consider the following time-varying optimization problem over a finite time duration $T\geq0$:
\begin{align}\label{eq:offline_problem}
  \min \sum_{t=0}^T f^t(x_t),\quad x_t \in \Omega,\quad t = 0,\ldots,T
\end{align}
where $\Omega\subseteq \mathbb{R}^p$ is a constraint set, $f^t:\mathbb{R}^p\to\mathbb{R}$ is the global objective function at time-step $t$. Moreover, we consider that the global objective function at each time-step $t$ is cast as the summation of $N$ local objective functions, given by
\begin{align}\label{eq:cost_function}
  f^t(x) = \sum_{i=1}^N f^t_i(x),
\end{align}
where $f^t_i:\mathbb{R}^p\to\mathbb{R}$ is a local objective function associated with agent $i$ at time-step $t$. 
At each time-step $t = 0,\ldots,T$, if the global objective function $f^t$ is known a priori, then one can directly obtain an optimizer $x^\star_t = \arg\min_{x \in \Omega} f^t(x)$ offline at each time-step $t=0,\ldots,T$. Hence, the optimal solution to problem~\eqref{eq:offline_problem} would be the sequence $\{x^\star_t\}_{t=0}^T$ with the optimal value being $\sum_{t=0}^Tf^t(x^\star_t)$. 
However, obtaining such solution sequence is infeasible due to the following settings: 1) (distributed setting) The global objective function is not directly accessible to any of the agents. 2) (online setting) At each time-step $t$, agent $i$ has no knowledge on the explicit expression of $f^t_i$ in the first place. After the decision $x^i_t$ is made by agent $i$ at time-step $t$, only the value of $f^t_i(x^i_t)$ with no information on its expression will then be revealed. In this case, at each time-step $t$, each agent can only generate $x^i_t$ based on the information from previous time-steps, and hence obtain a decision sequence $\{x^i_t\}_{t=0}^T$. Then, it would be interesting to analyze how close is the obtained decision sequence $\{x^i_t\}_{t=0}^T$ of each agent $i$ to the exact optimal solution sequence $\{x^\star_t\}_{t=0}^T$, by measuring the difference between their corresponding function values, which is known as dynamic regret. In particular, the dynamic regret due to agent $i \in \mathcal{V}$ is given by:
\begin{align*}
  \mathbf{R}_i(T) = \sum^T_{t=0}f^t(x^i_t)-\sum^T_{t=0}f^t(x^\star_t),
\end{align*}
where $x^i_t$ is agent $i$'s decision at time-step $t$. It should be noted that the dynamic regret $\mathbf{R}_i(T)$ depends on the problem regularities, where the most common one is the minimizer path length (\textit{i.e.}, the deviation of the consecutive optimal solution sequence), defined by:
\begin{align}\label{eq:def_opt_deviation}
  \omega_T = \sum^T_{t=0}\|x^\star_{t+1}-x^\star_t\|.
\end{align}
Thus, the aim of this work is to develop a distributed online optimization algorithm to solve problem~\eqref{eq:offline_problem} such that all agents' decision sequences are close to the optimal solution sequence.

We suppose the following assumption holds for the local objective function $f^t_i$:
\begin{Assumption}\label{assumption_local_f}
The local objective function $f^t_i:\mathbb{R}^p\to\mathbb{R}$, $i\in\mathcal{V}$ at time-step $t$ is convex and Lipschitz continuous in $\mathbb{R}^p$ with a constant $\hat{D}$, \textit{i.e.}, for $\forall x, x' \in \mathbb{R}^p$, we have $|f^t_i(x)-f^t_i(x')|\leq\hat{D}\|x-x'\|$. The constraint set $\Omega$ is convex and compact. 
\end{Assumption}

The interactions between agents are characterized by a directed graph $\mathcal{G} = \{\mathcal{V},\mathcal{E}\}$, where $\mathcal{V} = \{1, 2, \ldots, N\}$ denotes the set of agents, and $\mathcal{E} \subset \mathcal{V}\times\mathcal{V}$ is the set of edges. A directed edge ($i,j$) implies that agent $j$ can receive information from agent $i$. The set of in-neighbors (respectively, out-neighbors) of agent $i$ is denoted by $\mathcal{N}^{\text{in}}_i = \{j \in \mathcal{V} | (j,i)\in \mathcal{E}\}$ (respectively, $\mathcal{N}^{\text{out}}_i = \{j \in \mathcal{V} | (i,j)\in \mathcal{E}\}$), which contains agent $i$ itself. Regarding the network topology, we make the following standard assumption:
\begin{Assumption}\label{assumption_graph}
The digraph $\mathcal{G}$ is strongly connected.
\end{Assumption}

\subsection{Preliminaries}
To facilitate the gradient-free techniques, we formulate the following smoothed problem of \eqref{eq:offline_problem}
\begin{align*}
  \min f^t_\mu(x_t) = \sum_{i=1}^N f^t_{i,\mu^i}(x_t),\quad x_t \in \Omega,
\end{align*}
where $f^t_{i,\mu^i}(x_t)$ is a Gaussian approximation of $f^t_i(x_t)$
\begin{align*}
  f^t_{i,\mu^i}(x_t) = \frac1{\kappa}\int_{\mathbb{R}^p} f^t_i(x_t+\mu^i_t\xi^i_t)e^{-\frac12\|\xi^i_t\|^2}d\xi^i_t,
\end{align*}
with $\kappa = \int_{\mathbb{R}^p} e^{-\frac12\|\xi^i_t\|^2}d\xi^i_t = (2\pi)^{p/2}$ and $\mu^i_t \geq 0$ being a smoothing parameter sequence. Inspired by \cite{Nesterov2017}, the randomized gradient-free oracle of $f^t_i(x_t)$ at time-step $t$ is designed by
\begin{align}\label{prelim_oracle}
g^t_{\mu^i}(x_t) = \frac{f^t_i(x_t+\mu^i_t\xi^i_t)-f^t_i(x_t)}{\mu^i_t}\xi^i_t,
\end{align}
where $\xi^i_t \in \mathbb{R}^p$ is a normally distributed Gaussian vector. Let $\mathcal{F}_t$ denote the $\sigma$-field of all random variables up to time $t-1$.
Then, we can derive the following properties of the functions $g^t_{\mu^i}(x_t)$ and $f^t_{i,\mu^i}(x_t)$: 

\begin{Lemma}\label{lemma:property_f_mu}
(see \cite{Nesterov2017}) Let Assumption~\ref{assumption_local_f} hold. For $i \in \mathcal{V}$, functions $g^t_{\mu^i}(x_t)$ and $f^t_{i,\mu^i}(x_t)$ satisfy that
\begin{enumerate}
\item Function $f^t_{i,\mu^i}(x_t)$ is convex and differentiable with
\begin{align*}
f^t_i(x_t)\leq f^t_{i,\mu^i}(x_t)\leq f^t_i(x_t) + \sqrt{p}\hat{\mu}_t\hat{D},
\end{align*}
where $\hat{\mu}_t = \max_{i\in\mathcal{V}}\mu^i_t$.
\item The gradient $\nabla f^t_{i,\mu^i}(x_t)$ satisfies 
\begin{align*}
\nabla f^t_{i,\mu^i}(x_t) = \mathbf{E}[g^t_{\mu^i}(x_t)|\mathcal{F}_t],
\end{align*}
and is Lipschitz continuous, \textit{i.e.},
\begin{align*}
\|\nabla f^t_{i,\mu^i}(x_t) - \nabla f^t_{i,\mu^i}(y_t)\| \leq \hat{L}_t\|x_t - y_t\|,
\end{align*}
where $\hat{L}_t = \max_{i\in\mathcal{V}}L^i_t$, and $L^i_t=\frac{\sqrt{p}\hat{D}}{\mu^i_t}$.
\item The oracle $g^t_{\mu^i}(x_t)$ holds that
\begin{align*}
\mathbf{E}[\|g^t_{\mu^i}(x_t)\||\mathcal{F}_t] \leq \sqrt{\mathbf{E}[\|g^t_{\mu^i}(x_t)\|^2|\mathcal{F}_t]}\leq (p+4)\hat{D}.
\end{align*}
\end{enumerate}
\end{Lemma}

%%%%%%%%%%%%%%%%%%%%%%%%%%%%%%%%%%%%%%%%%%%%%%%%%%%%%%%%%%%%%%%%%%%%%
% \section{Main Results}\label{sec:distr_opt}

% In this section, the development of an online gradient-free distributed projected gradient descent method is presented, followed by the analysis on its convergence.

% \subsection{Online Gradient-Free DPGD Method} 

\section{Online Gradient-Free DPGD} \label{sec:distr_opt}

In this section, we propose an online gradient-free distributed projected gradient descent (DPGD) algorithm to solve the distributed online optimization problem, which is detailed as follows.

In this algorithm, every agent $i\in\mathcal{V}$ needs to maintain two variables, one state variable $x^i$ and one auxiliary variable $y^i$. At each time-step $t$, every agent $j\in\mathcal{V}$ broadcasts two pieces of information $x^j_t$ and $[\mathbf{W}_c]_{ij}y^j_t$ to its out-neighbor $i \in \mathcal{N}^\text{out}_j$. On receiving such information the in-neighbors, every agent $i$ updates its variables $x^i_{t+1}$ and $y^i_{t+1}$ referring to the following update law:
\begin{align}
x^i_{t+1} &= \Pi_{\Omega}\bigg[\sum_{j=1}^N [\mathbf{W}_r]_{ij}x^j_t + \delta y^i_t - \gamma_t g^t_{\mu^i}(x^i_t)\bigg],\nonumber\\
y^i_{t+1} &= \sum_{j=1}^N [\mathbf{W}_c]_{ij} y^j_t  - \sum_{j=1}^N [\mathbf{W}_r]_{ij}x^j_t +x^i_t- \delta y^i_t,\label{eq:algorithm}
\end{align}
where $g^t_{\mu^i}(x^i_t)$ is designed based on \eqref{prelim_oracle}, given by
\begin{align}
g^t_{\mu^i}(x^i_t) = \frac{f^t_i(x^i_t+\mu^i_t\xi^i_t)-f^t_i(x^i_t)}{\mu^i_t}\xi^i_t, \label{grad_oracle}
\end{align}
$\mathbf{W}_r,\mathbf{W}_c$ correspond to the weighting matrices of the digraph $\mathcal{G}$ with $\mathbf{W}_r$ being row-stochastic and $\mathbf{W}_c$ being column-stochastic, \textit{i.e.,} 
\begin{align*}
\sum_{j=1}^N[\mathbf{W}_r]_{ij} = 1, \forall i \in \mathcal{V}, \quad\sum_{i=1}^N[\mathbf{W}_c]_{ij} = 1, \forall j \in \mathcal{V},
\end{align*}
% $\sum_{j=1}^N[\mathbf{W}_r]_{ij} = 1$ for all $i \in \mathcal{V}$, and $\sum_{i=1}^N[\mathbf{W}_c]_{ij} = 1$ for all $j \in \mathcal{V}$, 
$\gamma_t > 0$ is a non-increasing step-size, and $\delta$ is a small positive number. 
% The proposed algorithm is detailed in Algorithm~\ref{algo:rgf_d_dgd}.

\section{Performance Analysis via Dynamic Regret}\label{sec:analysis}

This section presents the main results on the convergence analysis of the algorithm by deriving a bound on its dynamic regret. 

To facilitate the analysis, we write (\ref{eq:algorithm}) compactly as
\begin{align}\label{eq:d-dgd}
{\phi}^i_{t+1} = \sum_{j=1}^{2N}[\mathbf{W}]_{ij}{\phi}^j_t + \theta^i_t,
\end{align}
where ${\phi}^i_t = x^i_t$ for $i\in\{1,\ldots,N\}$ and ${\phi}^i_t = y^{i-N}_t$ for $i\in\{N+1,\ldots,2N\}$; $\theta^i_t = x^i_{t+1}-\sum_{j=1}^N [\mathbf{W}_r]_{ij}x^j_t - \delta y^i_t$ for $i\in\{1,\ldots,N\}$ and $\theta^i_t = \mathbf{0}_n$ for $i\in\{N+1,\ldots,2N\}$; and $\mathbf{W} = [\begin{smallmatrix} \mathbf{W}_r & \delta I \\ I-\mathbf{W}_r & \mathbf{W}_c - \delta I\end{smallmatrix}]$. Due to the compactness of $\Omega$, we define $\rho = \sup_{x\in\Omega}\|x\|$, which implies $\|{\phi}^i_t\| \leq \rho, i\in\{1,\ldots,N\}, t\geq0$.
Next, we introduce two results, one regarding the convergence of the augmented weighting matrix $\mathbf{W}$, and the other related to the boundeness of the augmented oracle $\theta^i_t$.

\begin{Lemma}\label{lemma:A_matrix}
(see \cite{Cai2012,Xi2016}) Let Assumption~\ref{assumption_graph} hold. The augmented weighting matrix $\mathbf{W}$ will converge at a geometric rate, if constant $\delta$ in $\mathbf{W}$ is selected within $(0,\hat{\delta})$, \text{i.e.,}
\begin{equation*}
\left\|\mathbf{W}^t - \begin{bmatrix} \frac{\mathbf{1}_N\mathbf{1}^T_N}N &\frac{\mathbf{1}_N\mathbf{1}^T_N}N \\ \mathbf{0} & \mathbf{0} \end{bmatrix}\right\|_\infty \leq {C} \lambda^t, \quad t\geq 1,
\end{equation*}
where $\hat{\delta} = (\frac{1-|\sigma_3|}{20+8N})^N$, $\sigma_3$ is the third largest eigenvalue of $\mathbf{W}$ by setting $\delta =0$. ${C} >0$ and $0 < \lambda < 1$ are some constants.
\end{Lemma}
\begin{Remark}
It is noted that the theoretical bound $\hat{\delta}$ is conservative to ensure the convergence of the weighting matrix $\mathbf{W}$ for any arbitrary strongly connected directed and fixed topologies. For special graph structures, such as symmetric digraphs, cyclic digraphs, \textit{etc}., the bound on $\delta$ can be less conservative. More details can be referred to the work \cite{Cai2012}.
\end{Remark}

For easy representation, we denote $\sum_{j=1}^{N}\|\theta^j_t\|$ by $\Theta_t$ in the rest of the paper.
\begin{Lemma}\label{lemma:bound_on_g}
With Assumptions~\ref{assumption_local_f} and \ref{assumption_graph}, let the constant $\delta$ be chosen such that $0<\delta <\min(\hat{\delta},\frac{1-\lambda}{2\sqrt{3}N{C}\lambda})$. Then, for any $T\geq0$, we have
\begin{enumerate}
\item $\begin{aligned}\sum_{t=0}^T\mathbf{E}[\Theta_t] \leq \frac{N(p+4)\hat{D}}{1-(1+2N\delta C)\lambda}\sum_{t=0}^T\gamma_t+H_1\end{aligned}$, 
\item $\begin{aligned}\sum_{t=0}^T\mathbf{E}[\Theta^2_t] \leq \frac{3N^2(p+4)^2\hat{D}^2}{1-(1+2\sqrt{3}N\delta C)\lambda}\sum_{t=0}^T\gamma^2_t +H_2\end{aligned}$,
\item $\begin{aligned}\sum_{t=0}^T\sqrt{\mathbf{E}[\Theta^2_t]} \leq \frac{\sqrt{3}N(p+4)\hat{D}}{1-(1+2\sqrt{3}N\delta C)\lambda}\sum_{t=0}^T\gamma_t+H_3\end{aligned}$,
\end{enumerate}
where $\hat{\delta}>0$, ${C}>0$ and $0<\lambda<1$ are the constants defined in Lemma~\ref{lemma:A_matrix}, $H_1$, $H_2$ and $H_3$ are bounded constants, and $\gamma_t>0$ is the step-size sequence.
\end{Lemma}
\begin{Proof}
See Appendix A for the proof.
\end{Proof}

We denote that $\bar{{\phi}}_t = \frac1N\sum_{i=1}^{2N}{\phi}^i_t$. Then, from \eqref{eq:d-dgd}, it yields that
\begin{align}
\bar{{\phi}}_{t+1} &= \frac1N\sum_{i=1}^{2N}\sum_{j=1}^{2N}[\mathbf{W}]_{ij}{\phi}^j_t + \frac1N\sum_{i=1}^{2N}\theta^i_t=\frac1N\sum_{j=1}^{2N}{\phi}^j_t+ \frac1N\sum_{i=1}^{2N}\theta^i_t=\bar{{\phi}}_t+\frac1N\sum_{i=1}^{2N}\theta^i_t, \label{eq:bar_z_relation}
\end{align}
where the column-stochastic property of $\mathbf{W}$ has been applied.
Then we provide bounds for $\mathbf{E}[\|{\phi}^i_t - \bar{{\phi}}_t\|]$ and $\mathbf{E}[\|{\phi}^i_t-\bar{{\phi}}_t\|^2]$, $1\leq i \leq N$.

\begin{Lemma}\label{lemma:consensus}
Suppose Assumptions~\ref{assumption_local_f} and \ref{assumption_graph} hold. Let the constant $\delta$ be chosen such that $0<\delta < \min(\hat{\delta},\frac{1-\lambda}{2\sqrt{3}N{C}\lambda})$. Let $\{{\phi}^i_t\}_{t\geq0}$ be the sequence generated by (\ref{eq:d-dgd}). Then, for $i = \{1,\ldots,N\}$, ${\phi}^i_t$ satisfies
\begin{align*}
\mathbf{E}[\|{\phi}^i_t-\bar{{\phi}}_t\|] &\leq 2N\rho\hat{{C}}\lambda^t+ \hat{{C}}\sum^{t}_{r=1}\lambda^{t-r}\mathbf{E}[\Theta_{r-1}],\\
\mathbf{E}[\|{\phi}^i_t-\bar{{\phi}}_t\|^2] &\leq 8N^2\rho^2\hat{{C}}^2\lambda^{2t}+\frac{2\hat{{C}}^2}{1-\lambda}\sum_{r=1}^t\lambda^{t-r}\mathbf{E}[\Theta^2_{r-1}],
\end{align*}
where $\hat{{C}}=\max\{{C},1\}>0$, $\hat{\delta}>0$, ${C}>0$ and $0<\lambda<1$ are constants defined in Lemma~\ref{lemma:A_matrix}.
\end{Lemma}
\begin{Proof}
Based on (\ref{eq:d-dgd}), we have that for $t \geq 1$ 
\begin{align}
{\phi}^i_t &= \sum_{j=1}^{2N}[\mathbf{W}^t]_{ij}{\phi}^j_0  + \sum_{r=1}^{t-1}\sum_{j=1}^{2N}[\mathbf{W}^{t-r}]_{ij}\theta^j_{r-1} + \theta^i_{t-1}.\label{eq:z}
\end{align}
It is noted that
\begin{equation}\label{eq:z_bar}
\bar{{\phi}}_t = \frac1N\sum_{j=1}^{2N}{\phi}^j_0 + \frac1N\sum_{r=1}^{t}\sum_{j=1}^{2N}\theta^j_{r-1}.
\end{equation}
Subtracting (\ref{eq:z_bar}) from (\ref{eq:z}) for $1 \leq i \leq N$, $t\geq 1$ gives
\begin{align*}
\|{\phi}^i_t-\bar{{\phi}}_t\| &\leq \sum_{j=1}^{2N}\bigg\|[\mathbf{W}^t]_{ij}-\frac1N\bigg\|\rho+ \sum_{r=1}^{t-1}\sum_{j=1}^{N}\bigg\|[\mathbf{W}^{t-r}]_{ij}-\frac1N\bigg\|\|\theta^j_{r-1}\|\\
&\quad+ \frac{N-1}N\|\theta^i_{t-1}\| + \frac1N\sum_{j\neq i}\|\theta^j_{t-1}\|.
\end{align*}
The last two terms satisfy
\begin{align*}
&\frac{N-1}N\|\theta^i_{t-1}\| + \frac1N\sum_{j\neq i}\|\theta^j_{t-1}\|\leq\frac{N-1}N\Theta_{t-1}+ \frac1N\Theta_{t-1}=\Theta_{t-1}.
\end{align*}
Then we obtain the following relation by Lemma~\ref{lemma:A_matrix} and defining $\hat{{C}}=\max\{{C},1\}$
\begin{align}
\|{\phi}^i_t-\bar{{\phi}}_t\| \leq2N\rho\hat{{C}}\lambda^t+\hat{{C}}\sum_{r=1}^{t}\lambda^{t-r}\Theta_{r-1}. \label{eq:consensus_error_relation}
\end{align}
Squaring both sides of \eqref{eq:consensus_error_relation} and applying Cauchy-Schwarz inequality
\begin{align*}
\bigg(\sum_{r=1}^t\lambda^{t-r}\Theta_{r-1}\bigg)^2 &\leq \bigg(\sum_{r=1}^t\lambda^{t-r}\bigg)\bigg(\sum_{r=1}^t\lambda^{t-r}\Theta^2_{r-1}\bigg)\leq\frac1{1-\lambda}\sum_{r=1}^t\lambda^{t-r}\Theta^2_{r-1},
\end{align*}
we obtain
\begin{align}
\|{\phi}^i_t-\bar{{\phi}}_t\|^2 &\leq 8N^2\rho^2\hat{{C}}^2\lambda^{2t}+\frac{2\hat{{C}}^2}{1-\lambda}\sum_{r=1}^t\lambda^{t-r}\Theta^2_{r-1}. \label{eq:consensus_error_relation_square}
\end{align}
% \begin{align}
% \|{\phi}^i_t-\bar{{\phi}}_t\| &\leq 2N\rho{C}\lambda^t+ \Theta_{t-1}+{C}\sum_{r=1}^{t-1}\lambda^{t-r}\Theta_{r-1}\nonumber\\
% &\leq2N\rho\hat{{C}}\lambda^t+\hat{{C}}\sum_{r=1}^{t}\lambda^{t-r}\Theta_{r-1}, \label{eq:consensus_error_relation}\\
% \|{\phi}^i_t-\bar{{\phi}}_t\|^2 &\leq \bigg(2N\rho\hat{{C}}\lambda^t+\hat{{C}}\sum_{r=1}^{t}\lambda^{t-r}\Theta_{r-1}\bigg)^2\nonumber\\
% &\leq 8N^2\rho^2\hat{{C}}^2\lambda^{2t}+\frac{2\hat{{C}}^2}{1-\lambda}\sum_{r=1}^t\lambda^{t-r}\Theta^2_{r-1}, \label{eq:consensus_error_relation_square}
% \end{align}
% where we have applied Cauchy-Schwarz inequality
% \begin{align*}
% \bigg(\sum_{r=1}^t\lambda^{t-r}\Theta_{r-1}\bigg)^2 &\leq \bigg(\sum_{r=1}^t\lambda^{t-r}\bigg)\bigg(\sum_{r=1}^t\lambda^{t-r}\Theta^2_{r-1}\bigg)\leq\frac1{1-\lambda}\sum_{r=1}^t\lambda^{t-r}\Theta^2_{r-1}.
% \end{align*}
Taking the total expectation yields the desired results.
\end{Proof}

% \begin{Remark}
% Denoting $\lim_{t\to\infty}\gamma_t$ by $\tilde{\gamma}$, then we have
% \begin{align*}
% \lim_{t\to\infty}\sum_{r=1}^{t}\lambda^{t-r}\gamma_{r-1} = \frac{\tilde{\gamma}}{1-\lambda},
% \end{align*}
% which means $\limsup_{t\to\infty}\mathbf{E}[\|{\phi}^i_t-\bar{{\phi}}_t\|]\leq(G_1\hat{{C}}\tilde{\gamma})/(1-\lambda)$; namely,
% all agents ${\phi}^i_t, i \in \mathcal{V}$ approximately converge to $\bar{{\phi}}_t$ with the error gap proportional to $\tilde{\gamma}$. Exact convergence is attained with $\tilde{\gamma} = 0$.
% \end{Remark}

Now we are ready to establish a bound for the dynamic regret $\mathbf{R}_i(T)$.
\begin{Theorem}\label{theorem:dynamic_regret}
Suppose Assumptions~\ref{assumption_local_f} and \ref{assumption_graph} hold. Let the constant $\delta$ be chosen such that $\delta \leq \min(\hat{\delta},\frac{1-\lambda}{2\sqrt{3}N{C}\lambda})$. Let $\{{\phi}^i_t\}_{t\geq0}$ be the sequence generated by (\ref{eq:d-dgd}) with a non-increasing step-size sequence $\{\gamma_t\}_{t\geq0}$ and a non-increasing smoothing parameter sequence $\{\mu^i_t\}_{t\geq0}$. Then, for $i\in\mathcal{V}$ and $T>0$, the dynamic regret $\mathbf{R}_i(T)$ holds that
\begin{align*}
\mathbf{E}[\mathbf{R}_i(T)]&\leq \mathcal{K}_1\sum^T_{t=0}\hat{\mu}_t + \frac{\mathcal{K}_2\omega_T+\mathcal{K}_3}{\gamma_T}+\bigg(\frac{\mathcal{K}_4}{\gamma_T}+\mathcal{K}_5\hat{L}_T\bigg)\sum_{t=0}^T\gamma^2_t\\
&\quad+(\mathcal{K}_6\hat{L}_T+\mathcal{K}_7)\sum_{t=0}^T\gamma_t+\mathcal{K}_8\hat{L}_T+\mathcal{K}_9,
\end{align*}
where $\mathcal{K}_1 = N\sqrt{p}\hat{D}$, $\mathcal{K}_2 = 2N\rho$, $\mathcal{K}_3 = N\hat{V}+\frac{N\rho\hat{{C}}}{1-\lambda^2}+H_2(N\rho\hat{{C}}+\frac{\hat{C}}{1-\lambda})$, $\mathcal{K}_4 = N\rho\hat{{C}}+\frac{\hat{C}}{1-\lambda}$, $\mathcal{K}_5 = \frac{2N\hat{{C}}^2}{(1-\lambda)^2}$, $\mathcal{K}_6 = \frac{2N^2(p+4)\rho\hat{D}\hat{{C}}}{(1-\lambda)(1-(1+2N\delta C)\lambda)}$, $\mathcal{K}_7 = \frac{\sqrt{3}N^2(p+4)^2\hat{D}^2\hat{C}}{(1-\lambda)(1-(1+2\sqrt{3}N\delta C)\lambda)}+ \frac{N^2(p+4)\hat{D}^2\hat{{C}}}{(1-\lambda)(1-(1+2N\delta C)\lambda)}+\frac{(p+4)^2N\hat{D}^2}2$, $\mathcal{K}_8 = \frac{2N\rho\hat{{C}}H_1}{1-\lambda}+\frac{2N\hat{{C}}^2H_2}{(1-\lambda)^2}+\frac{4N^2\rho^2\hat{{C}}}{1-\lambda}+\frac{8N^3\rho^2\hat{{C}}^2}{1-\lambda^2}$, $\mathcal{K}_9 = \frac{2N^2(p+5)\rho\hat{D}\hat{{C}}}{1-\lambda}+\frac{N\hat{D}\hat{{C}}H_1}{1-\lambda}+\frac{N(p+4)\hat{D}\hat{C}H_3}{1-\lambda}$,
$\hat{V}>0$ is a constant, $\rho = \sup_{x\in\Omega}\|x\|$, $\omega_T$ is defined in \eqref{eq:def_opt_deviation}, $\hat{L}_t>0$ is a constant defined in Lemma~\ref{lemma:property_f_mu}, $\hat{{C}}=\max\{{C},1\}>0$, $\hat{\delta}>0$, $C>0$ and $0<\lambda<1$ are constants defined in Lemma~\ref{lemma:A_matrix}, $H_1>0$, $H_2>0$ and $H_3>0$ are constants defined in Lemma~\ref{lemma:bound_on_g}.
\end{Theorem}
\begin{Proof}
Define a positive scalar function $V_t$ as
\begin{align*}
V_t = \frac12\langle\bar{{\phi}}_t-x^\star_t,\bar{{\phi}}_t-x^\star_t\rangle,
\end{align*}
where $\bar{{\phi}}_t$ is given in \eqref{eq:bar_z_relation} and $x^\star_t$ is the optimal solution of \eqref{eq:offline_problem} at time-step $t$. Then, it can be obtained that
\begin{align}
\Delta V_t &= \nu_{t+1} - V_t\nonumber\\
&= -\frac12\|\bar{{\phi}}_{t+1}-\bar{{\phi}}_t\|^2+\frac12\langle x^\star_{t+1}+x^\star_t-2\bar{{\phi}}_{t+1},x^\star_{t+1}-x^\star_t\rangle\nonumber\\
&\quad+\langle\bar{{\phi}}_{t+1}-\bar{{\phi}}_t,\bar{{\phi}}_{t+1}-x^\star_t\rangle\nonumber\\
&\leq -\frac12\|\bar{{\phi}}_{t+1}-\bar{{\phi}}_t\|^2+2\rho\|x^\star_{t+1}-x^\star_t\|\nonumber\\
&\quad+\langle\bar{{\phi}}_{t+1}-\bar{{\phi}}_t,\bar{{\phi}}_{t+1}-x^\star_t\rangle,\label{eq:prod_z_bar_minus_z_bar_times_z_bar_minus_xstar}
\end{align}
where we have used $\rho = \sup_{x\in\Omega}\|x\|$. Then, we will bound the dynamic regret $\mathbf{R}_i(T)$ in the following two steps.

\textbf{Step 1}. Bound of function $\Delta V_t$:

According to \eqref{eq:bar_z_relation}, we can further expand the last term in \eqref{eq:prod_z_bar_minus_z_bar_times_z_bar_minus_xstar} as follows
\begin{subequations}\label{eq:prod_z_bar_minus_z_bar_times_z_bar_minus_xstar_2}
\begin{align}
\langle\bar{{\phi}}_{t+1}-\bar{{\phi}}_t,\bar{{\phi}}_{t+1}-x^\star_t\rangle
&=\bigg\langle\frac1N\sum^{N}_{i=1}\theta^i_t,\bar{{\phi}}_{t+1}-x^\star_t\bigg\rangle\nonumber\\
&=\frac1N\sum^{N}_{i=1}\langle \theta^i_t,\bar{{\phi}}_{t+1}-{\phi}^i_{t+1}\rangle \label{eq:prod_g_times_z_bar_minus_zi}\\
&\quad+\frac1N\sum^{N}_{i=1}\langle \theta^i_t,{\phi}^i_{t+1}-x^\star_t\rangle \label{eq:prod_g_times_zi_minus_xstar}.
\end{align}
\end{subequations}

For \eqref{eq:prod_g_times_z_bar_minus_zi}, it yields that
\begin{align*}
\frac1N\sum^{N}_{i=1}\langle \theta^i_t,\bar{{\phi}}_{t+1}-{\phi}^i_{t+1}\rangle
&\leq\frac1N\sum^{N}_{i=1}\|\theta^i_t\|\|\bar{{\phi}}_{t+1}-{\phi}^i_{t+1}\|\\
&\leq\frac1N\Theta_t\bigg(2N\rho\hat{{C}}\lambda^{t+1}+\hat{{C}}\sum_{r=1}^{t+1}\lambda^{t-r+1}\Theta_{r-1}\bigg),
\end{align*}
where \eqref{eq:consensus_error_relation} has been substituted.
Taking the total expectation on both sides, and noting that $2\lambda^{t+1}\mathbf{E}[\Theta_t]\leq \lambda^{2(t+1)}+\mathbf{E}[\Theta^2_t]$, and $\mathbf{E}[\Theta_t\Theta_{r-1}]\leq\sqrt{\mathbf{E}[\Theta^2_t]\mathbf{E}[\Theta^2_{r-1}]}\leq\frac12(\mathbf{E}[\Theta^2_t]+\mathbf{E}[\Theta^2_{r-1}])$, we have
\begin{align}
\mathbf{E}\bigg[\frac1N\sum^{N}_{i=1}\langle \theta^i_t,\bar{{\phi}}_{t+1}-{\phi}^i_{t+1}\rangle\bigg]&\leq\rho\hat{{C}}\lambda^{2(t+1)}+\bigg(\rho\hat{{C}}+\frac{\hat{C}}{2N(1-\lambda)}\bigg)\mathbf{E}[\Theta^2_t]\nonumber\\
&\quad+\frac{\hat{C}}{2N}\sum_{r=1}^{t+1}\lambda^{t-r+1}\mathbf{E}[\Theta^2_{r-1}]. \label{eq:bound_on_prod_g_times_z_bar_minus_zi}
\end{align}

For \eqref{eq:prod_g_times_zi_minus_xstar}, it yields that
\begin{subequations}\label{eq:prod_g_times_zi_minus_xstar_2}
\begin{align}
\frac1N\sum^{N}_{i=1}\langle \theta^i_t,{\phi}^i_{t+1}-x^\star_t\rangle&=\frac1N\sum^{N}_{i=1}\langle \theta^i_t+\gamma_tg^t_{\mu^i}({\phi}^i_t),{\phi}^i_{t+1}-x^\star_t\rangle\label{eq:g_plus_alpha_g_mu_times_zi_minus_xstar}\\
&\quad+\frac1N\sum^{N}_{i=1}\langle\gamma_tg^t_{\mu^i}({\phi}^i_t),x^\star_t-{\phi}^i_{t+1}\rangle.\label{eq:alpha_g_mu_times_zi_minus_xstar}
\end{align}
\end{subequations}
For \eqref{eq:g_plus_alpha_g_mu_times_zi_minus_xstar}, applying the projection inequality \cite[Eq.~(1)]{Nedic2010} yields
\begin{align}
\langle \theta^i_t+\gamma_tg^t_{\mu^i}({\phi}^i_t),{\phi}^i_{t+1}-x^\star_t\rangle\leq0.\label{eq:bound_on_g_plus_alpha_g_mu_times_zi_minus_xstar}
\end{align}
For \eqref{eq:alpha_g_mu_times_zi_minus_xstar}, it is further expanded as
\begin{subequations}\label{eq:alpha_g_mu_times_zi_minus_xstar_2}
\begin{align}
&\frac1N\sum^{N}_{i=1}\langle\gamma_tg^t_{\mu^i}({\phi}^i_t),x^\star_t-{\phi}^i_{t+1}\rangle
=\frac1N\sum^{N}_{i=1}\langle\gamma_tg^t_{\mu^i}({\phi}^i_t),x^\star_t-\bar{{\phi}}_t\rangle\label{eq:alpha_g_mu_times_xstar_minus_z_bar}\\
&\quad\quad\quad\quad\quad\quad\quad\quad+\frac1N\sum^{N}_{i=1}\langle\gamma_tg^t_{\mu^i}({\phi}^i_t),\bar{{\phi}}_t-\bar{{\phi}}_{t+1}\rangle\label{eq:alpha_g_mu_times_z_bar_minus_z_bar}\\
&\quad\quad\quad\quad\quad\quad\quad\quad+\frac1N\sum^{N}_{i=1}\langle\gamma_tg^t_{\mu^i}({\phi}^i_t),\bar{{\phi}}_{t+1}-{\phi}^i_{t+1}\rangle.\label{eq:alpha_g_mu_times_z_bar_minus_zi}
\end{align}
\end{subequations}
For \eqref{eq:alpha_g_mu_times_xstar_minus_z_bar}, we have
\begin{align*}
\frac1N\sum^{N}_{i=1}\langle\gamma_tg^t_{\mu^i}({\phi}^i_t),x^\star_t-\bar{{\phi}}_t\rangle
&=\frac1N\sum^{N}_{i=1}\langle\gamma_t\big(g^t_{\mu^i}({\phi}^i_t)-\nabla f^t_{i,\mu^i}(\bar{{\phi}}_t)\big),(x^\star_t-{\phi}^i_t)+({\phi}^i_t-\bar{{\phi}}_t)\rangle\\
&\quad+\frac1N\sum^{N}_{i=1}\langle\gamma_t\nabla f^t_{i,\mu^i}(\bar{{\phi}}_t),x^\star_t-\bar{{\phi}}_t\rangle.
\end{align*}
Taking the conditional expectation on $\mathcal{F}_t$ and due to Lemma~\ref{lemma:property_f_mu}-2), we obtain
\begin{align*}
\mathbf{E}\bigg[\frac1N\sum^{N}_{i=1}\langle\gamma_tg^t_{\mu^i}({\phi}^i_t),x^\star_t-\bar{{\phi}}_t\rangle\bigg|\mathcal{F}_t\bigg]&\leq\frac1N\sum^{N}_{i=1}\gamma_t\|\nabla f^t_{i,\mu^i}({\phi}^i_t)-\nabla f^t_{i,\mu^i}(\bar{{\phi}}_t)\|(\|x^\star_t-{\phi}^i_t\|+\|{\phi}^i_t-\bar{{\phi}}_t\|)\\
&\quad+\frac1N\sum^{N}_{i=1}\langle\gamma_t\nabla f^t_{i,\mu^i}(\bar{{\phi}}_t),x^\star_t-\bar{{\phi}}_t\rangle.
\end{align*}
Noting that $\|x^\star_t-{\phi}^i_t\|\leq 2\rho$ and Lemma~\ref{lemma:property_f_mu}-2) on $\|\nabla f^t_{i,\mu^i}({\phi}^i_t)-\nabla f^t_{i,\mu^i}(\bar{{\phi}}_t)\|$, it gives that
\begin{align*}
\mathbf{E}\bigg[\frac1N\sum^{N}_{i=1}\langle\gamma_tg^t_{\mu^i}({\phi}^i_t),x^\star_t-\bar{{\phi}}_t\rangle\bigg|\mathcal{F}_t\bigg]&\leq\frac{2\rho\hat{L}_t\gamma_t}N\sum^{N}_{i=1}\|{\phi}^i_t-\bar{{\phi}}_t\|+\frac{\hat{L}_t\gamma_t}N\sum^{N}_{i=1}\|{\phi}^i_t-\bar{{\phi}}_t\|^2\\
&\quad+\frac{\gamma_t}N\big(f^t_{\mu}(x^\star_t)-f^t(\bar{{\phi}}_t)\big),
\end{align*}
where convexity of $f^t_{i,\mu^i}(x)$ and Lemma~\ref{lemma:property_f_mu}-1) have been used. Taking the total expectation and applying Lemma~\ref{lemma:consensus}, we have
\begin{align}
&\mathbf{E}\bigg[\frac1N\sum^{N}_{i=1}\langle\gamma_tg^t_{\mu^i}({\phi}^i_t),x^\star_t-\bar{{\phi}}_t\rangle\bigg]\leq\frac{2\rho\hat{L}_t\gamma_t}N\sum^{N}_{i=1}\bigg(2N\rho\hat{{C}}\lambda^t+ \hat{{C}}\sum^{t}_{r=1}\lambda^{t-r}\mathbf{E}[\Theta_{r-1}]\bigg)\nonumber\\
&\quad+\frac{\hat{L}_t\gamma_t}N\sum^{N}_{i=1}\bigg(8N^2\rho^2\hat{{C}}^2\lambda^{2t}+\frac{2\hat{{C}}^2}{1-\lambda}\sum_{r=1}^t\lambda^{t-r}\mathbf{E}[\Theta^2_{r-1}]\bigg)+\frac{\gamma_t}N\big(f^t_{\mu}(x^\star_t)-f^t(\bar{{\phi}}_t)\big)\nonumber\\
&\leq4N\rho^2\hat{{C}}\hat{L}_t\gamma_t\lambda^t+8N^2\rho^2\hat{{C}}^2\hat{L}_t\gamma_t\lambda^{2t}+2\rho\hat{{C}}\hat{L}_t\gamma_t\sum^t_{r=1}\lambda^{t-r}\mathbf{E}[\Theta_{r-1}]\nonumber\\
&\quad+\frac{2\hat{{C}}^2\hat{L}_t\gamma_t}{1-\lambda}\sum_{r=1}^t\lambda^{t-r}\mathbf{E}[\Theta^2_{r-1}]+\frac{\gamma_t}N\big(f^t_{\mu}(x^\star_t)-\mathbf{E}[f^t(\bar{{\phi}}_t)]\big).\label{eq:bound_on_alpha_g_mu_times_xstar_minus_z_bar}
\end{align}
For \eqref{eq:alpha_g_mu_times_z_bar_minus_z_bar}, we have
\begin{align*}
\frac1N\sum^{N}_{i=1}\langle\gamma_tg^t_{\mu^i}({\phi}^i_t),\bar{{\phi}}_t-\bar{{\phi}}_{t+1}\rangle
&\leq\frac12\bigg\|\frac1N\sum^{N}_{i=1}\gamma_tg^t_{\mu^i}({\phi}^i_t)\bigg\|^2+\frac12\|\bar{{\phi}}_t-\bar{{\phi}}_{t+1}\|^2\\
&\leq\frac1{2N}\gamma^2_t\sum^{N}_{i=1}\|g^t_{\mu^i}({\phi}^i_t)\|^2+\frac12\|\bar{{\phi}}_t-\bar{{\phi}}_{t+1}\|^2.
\end{align*}
Taking the total expectation and applying Lemma~\ref{lemma:property_f_mu}-3) on $\mathbf{E}[g^t_{\mu^i}({\phi}^i_t)\|^2]$ yields
\begin{align}
\mathbf{E}\bigg[\frac1N\sum^{N}_{i=1}\langle\gamma_tg^t_{\mu^i}({\phi}^i_t),\bar{{\phi}}_t-\bar{{\phi}}_{t+1}\rangle\bigg]&\leq\frac{(p+4)^2\hat{D}^2}2\gamma^2_t+\frac12\mathbf{E}[\|\bar{{\phi}}_t-\bar{{\phi}}_{t+1}\|^2].\label{eq:bound_on_alpha_g_mu_times_z_bar_minus_z_bar}
\end{align}
For \eqref{eq:alpha_g_mu_times_z_bar_minus_zi}, we have
\begin{align*}
&\frac1N\sum^{N}_{i=1}\langle\gamma_tg^t_{\mu^i}({\phi}^i_t),\bar{{\phi}}_{t+1}-{\phi}^i_{t+1}\rangle\leq\frac1N\sum^{N}_{i=1}\gamma_t\|g^t_{\mu^i}({\phi}^i_t)\|\|\bar{{\phi}}_{t+1}-{\phi}^i_{t+1}\|\\
&\quad\quad\quad\leq\frac1N\sum^{N}_{i=1}\gamma_t\|g^t_{\mu^i}({\phi}^i_t)\|\bigg(2N\rho\hat{{C}}\lambda^{t+1}+\hat{{C}}\sum_{r=1}^{t+1}\lambda^{t-r+1}\Theta_{r-1}\bigg),
\end{align*}
where \eqref{eq:consensus_error_relation} has been substituted.
Taking the total expectation on both sides yields
\begin{align}
&\mathbf{E}\bigg[\frac1N\sum^{N}_{i=1}\langle\gamma_tg^t_{\mu^i}({\phi}^i_t),\bar{{\phi}}_{t+1}-{\phi}^i_{t+1}\rangle\bigg]\leq2\rho\hat{C}\gamma_t\lambda^{t+1}\sum^{N}_{i=1}\mathbf{E}[\|g^t_{\mu^i}({\phi}^i_t)\|]\nonumber\\
&\quad\quad+\frac{\hat{C}}N\gamma_t\sum_{r=1}^{t+1}\lambda^{t-r+1}\sum^{N}_{i=1}\sqrt{\mathbf{E}[\|g^t_{\mu^i}({\phi}^i_t)\|^2]\mathbf{E}[\Theta^2_{r-1}]}\nonumber\\
&\quad\leq2N(p+4)\rho\hat{D}\hat{C}\gamma_t\lambda^{t+1}+ (p+4)\hat{D}\hat{C}\gamma_t\sum_{r=1}^{t+1}\lambda^{t-r+1}\sqrt{\mathbf{E}[\Theta^2_{r-1}]},\label{eq:bound_on_alpha_g_mu_times_z_bar_minus_zi}
\end{align}
where Lemma~\ref{lemma:property_f_mu}-3) and Cauchy-Schwarz inequality have been applied.
Taking the total expectation for \eqref{eq:prod_g_times_zi_minus_xstar_2}, and combining the results of \eqref{eq:bound_on_g_plus_alpha_g_mu_times_zi_minus_xstar} and \eqref{eq:alpha_g_mu_times_zi_minus_xstar_2} with substitutions of \eqref{eq:bound_on_alpha_g_mu_times_xstar_minus_z_bar}, \eqref{eq:bound_on_alpha_g_mu_times_z_bar_minus_z_bar} and \eqref{eq:bound_on_alpha_g_mu_times_z_bar_minus_zi}, we have
\begin{align}
&\mathbf{E}\bigg[\frac1N\sum^{N}_{i=1}\langle \theta^i_t,{\phi}^i_{t+1}-x^\star_t\rangle\bigg]\leq4N\rho^2\hat{{C}}\hat{L}_t\gamma_t\lambda^t+8N^2\rho^2\hat{{C}}^2\hat{L}_t\gamma_t\lambda^{2t}\nonumber\\
&\quad\quad+2\rho\hat{{C}}\hat{L}_t\gamma_t\sum^t_{r=1}\lambda^{t-r}\mathbf{E}[\Theta_{r-1}]+\frac{2\hat{{C}}^2\hat{L}_t\gamma_t}{1-\lambda}\sum_{r=1}^t\lambda^{t-r}\mathbf{E}[\Theta^2_{r-1}]\nonumber\\
&\quad\quad+\frac{\gamma_t}N(f^t_{\mu}(x^\star_t)-\mathbf{E}[f^t(\bar{{\phi}}_t)])+\frac{(p+4)^2\hat{D}^2}2\gamma^2_t+\frac12\mathbf{E}[\|\bar{{\phi}}_t-\bar{{\phi}}_{t+1}\|^2]\nonumber\\
&\quad\quad+2N(p+4)\rho\hat{D}\hat{{C}}\gamma_t\lambda^{t+1}+(p+4)\hat{D}\hat{C}\gamma_t\sum_{r=1}^{t+1}\lambda^{t-r+1}\sqrt{\mathbf{E}[\Theta^2_{r-1}]}\nonumber\\
&\quad\leq2N(p+4)\rho\hat{D}\hat{{C}}\gamma_t\lambda^{t+1}+4N\rho^2\hat{{C}}\hat{L}_t\gamma_t\lambda^t+8N^2\rho^2\hat{{C}}^2\hat{L}_t\gamma_t\lambda^{2t}\nonumber\\
&\quad\quad+(p+4)\hat{D}\hat{C}\gamma_t\sum_{r=1}^{t+1}\lambda^{t-r+1}\sqrt{\mathbf{E}[\Theta^2_{r-1}]}+2\rho\hat{{C}}\hat{L}_t\gamma_t\sum^t_{r=1}\lambda^{t-r}\mathbf{E}[\Theta_{r-1}]\nonumber\\
&\quad\quad+\frac{2\hat{{C}}^2\hat{L}_t\gamma_t}{1-\lambda}\sum_{r=1}^t\lambda^{t-r}\mathbf{E}[\Theta^2_{r-1}]+\frac{(p+4)^2\hat{D}^2}2\gamma^2_t+\frac{\gamma_t}N(f^t_{\mu}(x^\star_t)-\mathbf{E}[f^t(\bar{{\phi}}_t)])\nonumber\\
&\quad\quad+\frac12\mathbf{E}[\|\bar{{\phi}}_t-\bar{{\phi}}_{t+1}\|^2].\label{eq:bound_on_prod_g_times_zi_minus_xstar}
\end{align}
Taking the total expectation for \eqref{eq:prod_z_bar_minus_z_bar_times_z_bar_minus_xstar_2}, and substituting \eqref{eq:bound_on_prod_g_times_z_bar_minus_zi} and \eqref{eq:bound_on_prod_g_times_zi_minus_xstar} into \eqref{eq:prod_g_times_z_bar_minus_zi} and \eqref{eq:prod_g_times_zi_minus_xstar}, we have
\begin{align}
&\mathbf{E}[\langle\bar{{\phi}}_{t+1}-\bar{{\phi}}_t,\bar{{\phi}}_{t+1}-x^\star_t\rangle]
\leq\rho\hat{{C}}\lambda^{2(t+1)}+2N(p+4)\rho\hat{D}\hat{{C}}\gamma_t\lambda^{t+1}+4N\rho^2\hat{{C}}\hat{L}_t\gamma_t\lambda^t\nonumber\\
&\quad+8N^2\rho^2\hat{{C}}^2\hat{L}_t\gamma_t\lambda^{2t}+ \bigg(\rho\hat{{C}}+\frac{\hat{C}}{2N(1-\lambda)}\bigg)\mathbf{E}[\Theta^2_t]+(p+4)\hat{D}\hat{C}\gamma_t\sum_{r=1}^{t+1}\lambda^{t-r+1}\sqrt{\mathbf{E}[\Theta^2_{r-1}]}\nonumber\\
&\quad+2\rho\hat{{C}}\hat{L}_t\gamma_t\sum^t_{r=1}\lambda^{t-r}\mathbf{E}[\Theta_{r-1}]+\frac{2\hat{{C}}^2\hat{L}_t\gamma_t}{1-\lambda}\sum_{r=1}^t\lambda^{t-r}\mathbf{E}[\Theta^2_{r-1}]+\frac{\hat{C}}{2N}\sum_{r=1}^{t+1}\lambda^{t-r+1}\mathbf{E}[\Theta^2_{r-1}]\nonumber\\
&\quad+\frac{(p+4)^2\hat{D}^2}2\gamma^2_t+\frac{\gamma_t}N(f^t_{\mu}(x^\star_t)-\mathbf{E}[f^t(\bar{{\phi}}_t)])+\frac12\mathbf{E}[\|\bar{{\phi}}_t-\bar{{\phi}}_{t+1}\|^2]. \label{eq:bound_on_prod_z_bar_minus_z_bar_times_z_bar_minus_xstar}
\end{align}
Taking the total expectation for \eqref{eq:prod_z_bar_minus_z_bar_times_z_bar_minus_xstar}, and substituting \eqref{eq:bound_on_prod_z_bar_minus_z_bar_times_z_bar_minus_xstar} into it, we can obtain
\begin{align}
&\mathbf{E}[\Delta V_t]\leq 2\rho\|x^\star_{t+1}-x^\star_t\|+\rho\hat{{C}}\lambda^{2(t+1)}+2N(p+4)\rho\hat{D}\hat{{C}}\gamma_t\lambda^{t+1}+4N\rho^2\hat{{C}}\hat{L}_t\gamma_t\lambda^t\nonumber\\
&\quad+8N^2\rho^2\hat{{C}}^2\hat{L}_t\gamma_t\lambda^{2t}+ \bigg(\rho\hat{{C}}+\frac{\hat{C}}{2N(1-\lambda)}\bigg)\mathbf{E}[\Theta^2_t]+(p+4)\hat{D}\hat{C}\gamma_t\sum_{r=1}^{t+1}\lambda^{t-r+1}\sqrt{\mathbf{E}[\Theta^2_{r-1}]}\nonumber\\
&\quad+2\rho\hat{{C}}\hat{L}_t\gamma_t\sum^t_{r=1}\lambda^{t-r}\mathbf{E}[\Theta_{r-1}]+\frac{2\hat{{C}}^2\hat{L}_t\gamma_t}{1-\lambda}\sum_{r=1}^t\lambda^{t-r}\mathbf{E}[\Theta^2_{r-1}]+\frac{\hat{C}}{2N}\sum_{r=1}^{t+1}\lambda^{t-r+1}\mathbf{E}[\Theta^2_{r-1}]\nonumber\\
&\quad+\frac{(p+4)^2\hat{D}^2}2\gamma^2_t+\frac{\gamma_t}N(f^t_{\mu}(x^\star_t)-\mathbf{E}[f^t(\bar{{\phi}}_t)]). \label{eq:bound_of_EVT}
\end{align}
Due to the compactness of $\Omega$, we define an upper bound of $V_t = \frac12\|\bar{{\phi}}_t-x^\star_t\|^2$ by $\hat{V}$. Then
\begin{align*}
-\sum^T_{t=0}\frac{\mathbf{E}[\Delta V_t]}{\gamma_t} &= \sum^T_{t=0}\frac{\mathbf{E}[V_t] - \mathbf{E}[\nu_{t+1}]}{\gamma_t}\leq\frac{\mathbf{E}[\nu_0]}{\gamma_0}+\sum^T_{t=1}\bigg(\frac1{\gamma_t}-\frac1{\gamma_{t-1}}\bigg)\mathbf{E}[V_t]\\
&\leq \frac{\hat{V}}{\gamma_0}+\hat{V}\sum^T_{t=1}\bigg(\frac1{\gamma_t}-\frac1{\gamma_{t-1}}\bigg)=\frac{\hat{V}}{\gamma_T},
\end{align*}
where we have used the fact that $\gamma_t$ is positive and non-increasing. Dividing both sides of \eqref{eq:bound_of_EVT} by $\gamma_t$, summing up from $t = 0$ to $T$ and combining the above relation, we have
\begin{align*}
&\sum^T_{t=0}(\mathbf{E}[f^t(\bar{{\phi}}_t)]-f^t_{\mu}(x^\star_t))\leq \frac{N\hat{V}}{\gamma_T}+\frac{2N\rho\omega_T}{\gamma_T}+2N^2(p+4)\rho\hat{D}\hat{{C}}\sum^T_{t=0}\lambda^{t+1}+4N^2\rho^2\hat{{C}}\hat{L}_T\sum^T_{t=0}\lambda^t\\
&\quad+\bigg(\frac{N\rho\hat{{C}}}{\gamma_T}+8N^3\rho^2\hat{{C}}^2\hat{L}_T\bigg)\sum^T_{t=0}\lambda^{2t}+ \bigg(\frac{N\rho\hat{{C}}}{\gamma_T}+\frac{\hat{C}}{(1-\lambda)\gamma_T}+\frac{2N\hat{{C}}^2\hat{L}_T}{(1-\lambda)^2}\bigg)\sum^T_{t=0}\mathbf{E}[\Theta^2_t]\\
&\quad+\frac{N(p+4)\hat{D}\hat{C}}{1-\lambda}\sum^T_{t=0}\sqrt{\mathbf{E}[\Theta^2_t]}+\frac{2N\rho\hat{{C}}\hat{L}_T}{1-\lambda}\sum^T_{t=0}\mathbf{E}[\Theta_t]+\frac{N(p+4)^2\hat{D}^2}2\sum^T_{t=0}\gamma_t,
\end{align*}
where we have used that $\gamma_t$ is non-increasing and $\hat{L}_t$ is non-decreasing as $\mu^i_t$ is non-increasing, and applied the following inequality ($A_t\geq0$)
\begin{align*}
&\sum^T_{t=0}\sum^{t+1}_{r=1}\lambda^{t-r+1}A_{r-1}\leq \frac1{1-\lambda}\sum^T_{t=0}A_t.
\end{align*}

\textbf{Step 2}. Bound of dynamic regret $\mathbf{R}_i(T)$:

% Denoting $\sum_{t=0}^T(B_1+B_2)$ by $B_T$ and using the result of $f^t_{\mu}(x)\leq f^t(x) + \sqrt{p}N\hat{\mu}_t\hat{D}$ from Lemma~\ref{lemma:property_f_mu}-1), we can obtain
Noting that $f^t_\mu(x)\leq f^t(x)+\sqrt{p}N\hat{\mu}_t\hat{D}$ from Lemma~\ref{lemma:property_f_mu}-1), we can obtain
\begin{align}
&\sum^T_{t=0}(\mathbf{E}[f^t(\bar{{\phi}}_t)]-f^t(x^\star_t)\leq N\sqrt{p}\hat{D}\sum^T_{t=0}\hat{\mu}_t+\frac{N\hat{V}+2N\rho\omega_T}{\gamma_T}+2N^2(p+4)\rho\hat{D}\hat{{C}}\sum^T_{t=0}\lambda^{t+1}\nonumber\\
&\quad+4N^2\rho^2\hat{{C}}\hat{L}_T\sum^T_{t=0}\lambda^t+\bigg(\frac{N\rho\hat{{C}}}{\gamma_T}+8N^3\rho^2\hat{{C}}^2\hat{L}_T\bigg)\sum^T_{t=0}\lambda^{2t}\nonumber\\
&\quad+ \bigg(\frac{N\rho\hat{{C}}}{\gamma_T}+\frac{\hat{C}}{(1-\lambda)\gamma_T}+\frac{2N\hat{{C}}^2\hat{L}_T}{(1-\lambda)^2}\bigg)\sum^T_{t=0}\mathbf{E}[\Theta^2_t]+\frac{N(p+4)\hat{D}\hat{C}}{1-\lambda}\sum^T_{t=0}\sqrt{\mathbf{E}[\Theta^2_t]}\nonumber\\
&\quad+\frac{2N\rho\hat{{C}}\hat{L}_T}{1-\lambda}\sum^T_{t=0}\mathbf{E}[\Theta_t]+\frac{N(p+4)^2\hat{D}^2}2\sum^T_{t=0}\gamma_t.\label{eq:bound_on_f_z_bar_minus_f_x_star}
\end{align}
Considering the dynamic regret $\mathbf{R}_i(T)$
% \footnote{The dynamic regret considered in this work is in the sense of expectation due to the agent's decision sequence $x^i_t$ generated by the proposed algorithm, which involves random variables in the estimation of gradient.}
\begin{align*}
\mathbf{E}[\mathbf{R}_i(T)] &=\sum^T_{t=0}\Big(\mathbf{E}[f^t({\phi}^i_t)]-f^t(x^\star_t)\Big)=\sum^T_{t=0}\Big(\mathbf{E}[f^t({\phi}^i_t)-f^t(\bar{{\phi}}_t)]\Big)\\
&\quad\quad\quad+\sum^T_{t=0}\Big(\mathbf{E}[f^t(\bar{{\phi}}_t)]-f^t(x^\star_t)\Big)\\
&\leq N\hat{D}\sum^T_{t=0}\mathbf{E}[\|{\phi}^i_t-\bar{{\phi}}_t\|]+\sum^T_{t=0}\Big(\mathbf{E}[f^t(\bar{{\phi}}_t)]-f^t(x^\star_t)\Big)\\
&\leq 2N^2\rho\hat{D}\hat{{C}}\sum^T_{t=0}\lambda^t+ \frac{N\hat{D}\hat{{C}}}{1-\lambda}\sum^T_{t=0}\mathbf{E}[\Theta_t]+\sum^T_{t=0}\Big(\mathbf{E}[f^t(\bar{{\phi}}_t)]-f^t(x^\star_t)\Big),
\end{align*}
where the first inequality follows from Assumption~\ref{assumption_local_f}, and the second inequality is due to Lemma~\ref{lemma:consensus}. Substituting \eqref{eq:bound_on_f_z_bar_minus_f_x_star} to the above gives,
\begin{align*}
&\mathbf{E}[\mathbf{R}_i(T)]\leq N\sqrt{p}\hat{D}\sum^T_{t=0}\hat{\mu}_t+\frac{N\hat{V}+2N\rho\omega_T}{\gamma_T}+\frac{N\rho\hat{{C}}}{(1-\lambda^2)\gamma_T}+\frac{2N^2(p+5)\rho\hat{D}\hat{{C}}}{1-\lambda}\\
&\quad+\frac{4N^2\rho^2\hat{{C}}\hat{L}_T}{1-\lambda}+\frac{8N^3\rho^2\hat{{C}}^2\hat{L}_T}{1-\lambda^2}+ \bigg(\frac{N\rho\hat{{C}}}{\gamma_T}+\frac{\hat{C}}{(1-\lambda)\gamma_T}+\frac{2N\hat{{C}}^2\hat{L}_T}{(1-\lambda)^2}\bigg)\sum^T_{t=0}\mathbf{E}[\Theta^2_t]\\
&\quad+\frac{N(p+4)\hat{D}\hat{C}}{1-\lambda}\sum^T_{t=0}\sqrt{\mathbf{E}[\Theta^2_t]}+\bigg(\frac{2N\rho\hat{{C}}\hat{L}_T}{1-\lambda}+ \frac{N\hat{D}\hat{{C}}}{1-\lambda}\bigg)\sum^T_{t=0}\mathbf{E}[\Theta_t]+\frac{N(p+4)^2\hat{D}^2}2\sum^T_{t=0}\gamma_t.
\end{align*}
Invoking Lemma~\ref{lemma:bound_on_g}, we have
\begin{align*}
&\mathbf{E}[\mathbf{R}_i(T)]\leq N\sqrt{p}\hat{D}\sum^T_{t=0}\hat{\mu}_t+\frac{N\hat{V}+2N\rho\omega_T}{\gamma_T}+\frac{N\rho\hat{{C}}}{(1-\lambda^2)\gamma_T}+\frac{2N^2(p+5)\rho\hat{D}\hat{{C}}}{1-\lambda}\\
&\quad+ \bigg(\frac{N\rho\hat{{C}}}{\gamma_T}+\frac{\hat{C}}{(1-\lambda)\gamma_T}+\frac{2N\hat{{C}}^2\hat{L}_T}{(1-\lambda)^2}\bigg)\bigg(\frac{3N^2(p+4)^2\hat{D}^2}{1-(1+2\sqrt{3}N\delta C)\lambda}\sum_{t=0}^T\gamma^2_t +H_2\bigg)\\
&\quad+\frac{N(p+4)\hat{D}\hat{C}}{1-\lambda}\bigg(\frac{\sqrt{3}N(p+4)\hat{D}}{1-(1+2\sqrt{3}N\delta C)\lambda}\sum_{t=0}^T\gamma_t+H_3\bigg)+\frac{4N^2\rho^2\hat{{C}}\hat{L}_T}{1-\lambda}+\frac{8N^3\rho^2\hat{{C}}^2\hat{L}_T}{1-\lambda^2}\\
&\quad+\bigg(\frac{2N\rho\hat{{C}}\hat{L}_T}{1-\lambda}+ \frac{N\hat{D}\hat{{C}}}{1-\lambda}\bigg)\bigg(\frac{N(p+4)\hat{D}}{1-(1+2N\delta C)\lambda}\sum_{t=0}^T\gamma_t+H_1\bigg)+\frac{N(p+4)^2\hat{D}^2}2\sum^T_{t=0}\gamma_t,
\end{align*}
which leads to the desired result.
\end{Proof}

\begin{Corollary}\label{corollary:dynamic_regret}
Suppose Assumptions~\ref{assumption_local_f} and \ref{assumption_graph} hold. For $i\in\mathcal{V}$ and $t = 0, \ldots, T$, let the step-size be $\gamma_t =\gamma = \frac{\gamma_0}{(T+1)^\alpha}$ and the smoothing parameter be $\mu^i_t=\mu^i = \frac{\mu_0}{(T+1)^\beta}$ with $\alpha,\beta \in (0,1)$. Then, the dynamic regret holds that
\begin{align*}
\mathbf{E}[\mathbf{R}_i(T)] &\leq \mathcal{O}(\max\{\omega_TT^\alpha,T^{\max\{1-\beta,1-\alpha+\beta\}}\}).
\end{align*}
Furthermore, if $\omega_T$ can be known in advance, by choosing $\gamma_0 = \omega_T^{\frac23}$, $\mu_0 = \omega_T^{\frac13}$, $\alpha = \frac23$ and $\beta = \frac13$, the dynamic regret can be improved as
\begin{align*}
\mathbf{E}[\mathbf{R}_i(T)] &\leq \mathcal{O}(\omega_T^{\frac13}T^{\frac23}).
\end{align*}
\end{Corollary}
\begin{Proof}
With $\gamma = \frac{\gamma_0}{(T+1)^\alpha}$ and $\mu^i = \frac{\mu_0}{(T+1)^\beta}$, it follows from Theorem~\ref{theorem:dynamic_regret} that
\begin{align*}
\mathbf{E}[\mathbf{R}_i(T)]&\leq \mathcal{K}_1\mu_0(T+1)^{1-\beta} + \frac{\mathcal{K}_2\omega_T+\mathcal{K}_3}{\gamma_0}(T+1)^\alpha+\mathcal{K}_4\gamma_0(T+1)^{1-\alpha}\\
&\quad+\frac{\mathcal{K}_5\sqrt{p}\hat{D}\gamma_0^2}{\mu_0}(T+1)^{1-2\alpha+\beta}+\frac{\mathcal{K}_6\sqrt{p}\hat{D}\gamma_0}{\mu_0}(T+1)^{1-\alpha+\beta}+\mathcal{K}_7\gamma_0(T+1)^{1-\alpha}\\
&\quad+\frac{\mathcal{K}_8\sqrt{p}\hat{D}}{\mu_0}(T+1)^\beta+\mathcal{K}_9.
\end{align*}
Noting that $\mathcal{O}(T^{1-\alpha})\leq\mathcal{O}(T^{1-\alpha+\beta})$, $\mathcal{O}(T^{1-2\alpha+\beta})\leq\mathcal{O}(T^{1-\alpha+\beta})$, $\mathcal{O}(T^\beta)\leq\mathcal{O}(T^{1-\alpha+\beta})$, $\mathcal{O}(T^\alpha)\leq\mathcal{O}(\omega_TT^\alpha)$, the bound of the dynamic regret can be written as 
\begin{align}
\mathbf{E}[\mathbf{R}_i(T)] \leq \mathcal{O}\bigg(\max\bigg\{\frac{\omega_TT^\alpha}{\gamma_0},\mu_0T^{1-\beta},\frac{\gamma_0}{\mu_0}T^{1-\alpha+\beta}\bigg\}\bigg). \label{eq:dynamic_regret_step_smooth_inequality}
\end{align}
which yields the first result.

Suppose for fixed $T>0$, $\omega_T$ can be known in advance. We may set $\gamma_0 = \omega_T^{\frac23}$, $\mu_0 = \omega_T^{\frac13}$, $\alpha = \frac23$ and $\beta = \frac13$. Then \eqref{eq:dynamic_regret_step_smooth_inequality} directly leads to the second result.
\end{Proof}

\begin{Remark}\label{remark:dynamic_regret}
To investigate the dependence of the dynamic regret on the number of agents $N$ and the problem dimension $p$, by noting that the coefficients $\mathcal{K}_1\leq\mathcal{O}(N\sqrt{p})$, $\mathcal{K}_2\leq\mathcal{O}(N)$, $\mathcal{K}_6\leq\mathcal{O}(N^2p)$, the dynamic regret in \eqref{eq:dynamic_regret_step_smooth_inequality} should be
\begin{align*}
\mathbf{E}[\mathbf{R}_i(T)] \leq \mathcal{O}\bigg(\max\bigg\{\frac{N}{\gamma_0}\omega_TT^\alpha,N\sqrt{p}\mu_0T^{1-\beta},\frac{N^2p\sqrt{p}\gamma_0}{\mu_0}T^{1-\alpha+\beta}\bigg\}\bigg). 
\end{align*}
In the case where $\omega_T$ can be known in advance, we may set $\gamma_0 = N^{-\frac13}p^{-\frac23}\omega_T^{\frac23}$, $\mu_0 = N^{\frac13}p^{\frac16}\omega_T^{\frac13}$, $\alpha = \frac23$ and $\beta = \frac13$. Then the dynamic regret becomes $\mathbf{E}[\mathbf{R}_i(T)] \leq \mathcal{O}(N^{\frac43}p^{\frac23}\omega_T^{\frac13}T^{\frac23})$.
\end{Remark}

%%%%%%%%%%%%%%%%%%%%%%%%%%%%%%%%%%%%%%%%%%%%%%%%%%%%%%%%%%%%%%%%%%%%%%%%%%%%%%%%
\section{Numerical Simulation}\label{sec:simulation}

\begin{figure}[!t]
\centering
\includegraphics[width=3in]{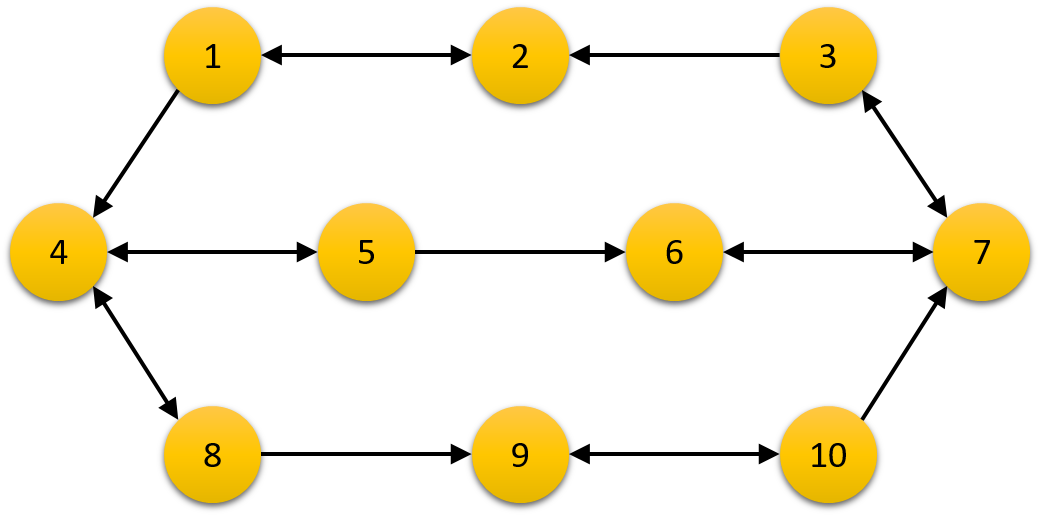}  % The printed column width is 8.4 cm.
\caption{Topology of the multi-agent system}
\label{fig:network.png}
\end{figure}

\begin{figure}[!t]
\centering
\includegraphics[width=4in]{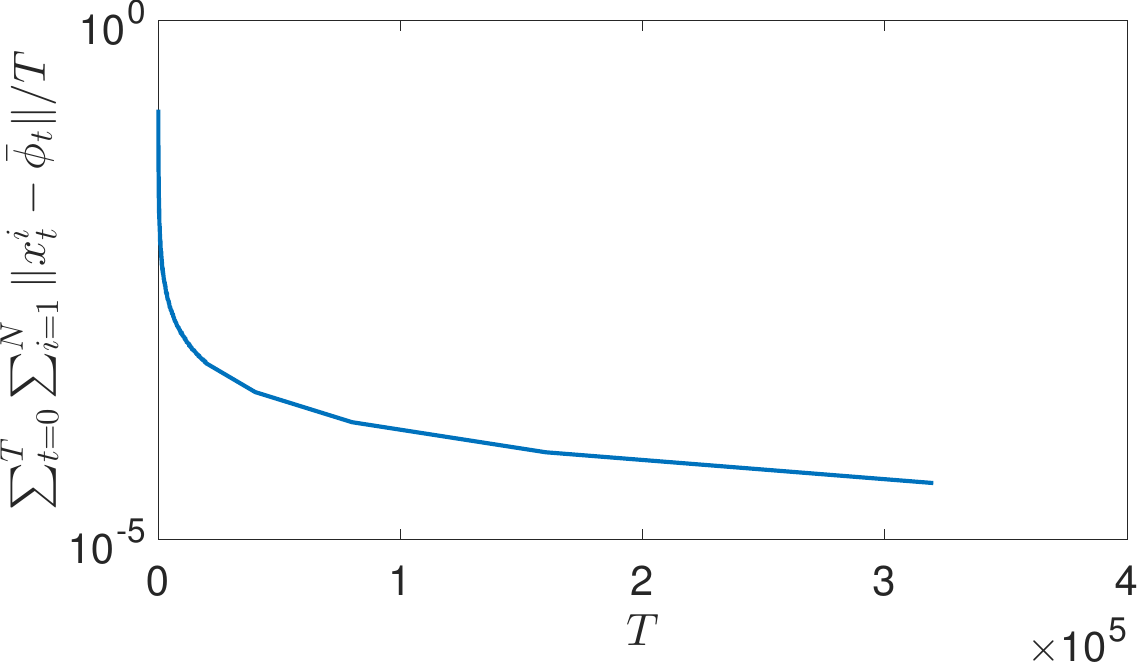}  % The printed column width is 8.4 cm.
\caption{Trajectories of $\sum_{t=0}^T\sum_{i=1}^{10}\mathbf{E}[\|x^i_t-\bar{{\phi}}_t\|]/T$.}
\label{fig:accumulation_error.eps}
\end{figure}

\begin{figure}[!t]
\centering
\includegraphics[width=4in]{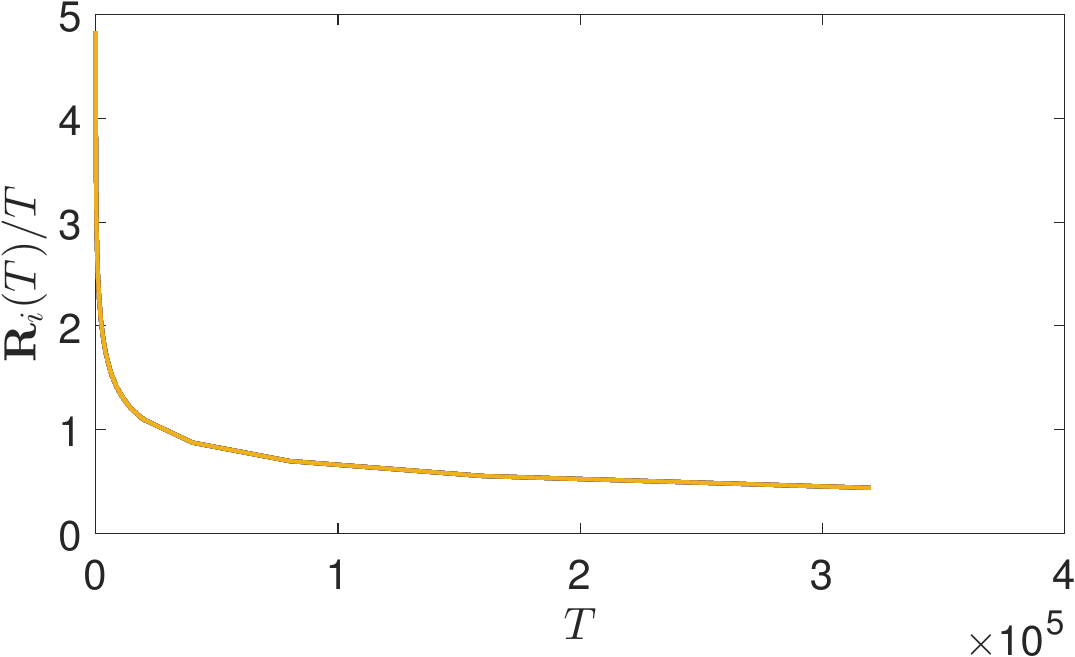}  % The printed column width is 8.4 cm.
\caption{Trajectories of $\mathbf{R}_i(T)/T$, $i=1,\ldots,10$}
\label{fig:average_dynamic_regret.eps}
\end{figure}

This section demonstrates the derived properties using a numerical example. In particular, we consider the following distributed online optimization problem with 10 agents over a strongly connected digraph as shown in Fig.~\ref{fig:network.png}:
\begin{equation*}
\min f^t(x) = \sum_{i=1}^{10}f^t_i(x), x\in\Omega,
\end{equation*}
where $f^t_i(x) = a_ix^2 - 2b_i(\frac{2\sin(0.008t)}t)x + c_i(\frac{2\sin(0.008t)}t)^2$ at time-step $t$, $\mathcal{X} = [-5,5]$, $a_i,b_i,c_i>0$ and $\sum_{i=1}^{10}a_i=\sum_{i=1}^{10}b_i=\sum_{i=1}^{10}c_i = 10$. Our proposed online gradient-free DPGD method will be used to solve this problem. For the weighting matrix $\mathbf{W}_r$ and $\mathbf{W}_c$, we let $[\mathbf{W}_r]_{ij} = 1/|\mathcal{N}^{\text{in}}_i|$ and $[\mathbf{W}_c]_{ij} = 1/|\mathcal{N}^{\text{out}}_j|$. We set $\delta = 0.1$, the smoothing parameter $\mu^i_t=\mu^i=1/(1+T)^{\frac13}$ and the step-size $\gamma_t =\gamma = 1/(1+T)^{\frac23}$ for $T = 1,\ldots, 3.2\times10^5$. 
% For any $t>0$, it is obvious that the optimal solution at time-step $t$ is $x^\star_t = \frac{2\sin(0.008t)}t$, which is time-varying. As $t$ increases, the deviation of the optimal solution $x^\star_t$ decreases gradually. We applied the proposed algorithm to solve the problem.
The trajectories of the time averaged consensus error ($\sum_{t=0}^T\sum_{i=1}^{10}\mathbf{E}[\|x^i_t-\bar{{\phi}}_t\|]/T$) and dynamic regret ($\mathbf{R}_i(T)/T$) were plotted in Figs.~\ref{fig:accumulation_error.eps} and \ref{fig:average_dynamic_regret.eps}, respectively. As can be seen, both the time averaged consensus error and dynamic regret gradually descend as expected, implying their sublinear growths.

% \begin{figure}[!t]
% \centering
% \includegraphics[width=3.4in]{convergence_of_x-eps-converted-to.pdf}  % The printed column width is 8.4 cm.
% \caption{Trajectories of $x^i_t$}
% \label{fig:convergence_of_x.eps}
% \end{figure}

\section{Conclusions}\label{sec:conclusion}
This work has considered a distributed online optimization problem. To solve the problem, we have developed an online gradient-free distributed projected gradient descent algorithm. A local randomized gradient-free oracle is built to estimate the gradient information in guiding the update of the decision variables. With some standard assumptions on graph connectivity and the objective function, we have quantified a bound of the dynamic regret for any agent as a function of the minimizer path length, step-size and smoothing parameter. Under appropriate selections of the step-size and smoothing parameter, we have proved that the dynamic regret is sublinear if the minimizer path length also grows sublinearly. Finally, numerical simulations have been conducted to verify the effectiveness of the algorithm.

% \begin{figure}
% \begin{center}
% \includegraphics[height=4cm]{jcaesar.eps}    % The printed column  
% \caption{Gaius Julius Caesar, 100--44 B.C.}  % width is 8.4 cm.
% \label{fig1}                                 % Size the figures 
% \end{center}                                 % accordingly.
% \end{figure}

% OR

%\begin{figure}
%\begin{center}
%\epsfig{file=jcaesar,width=7cm}
%\caption{Gaius Julius Caesar, 100--44 B.C.}
%\label{fig1}
%\end{center}
%\end{figure}

% \begin{ack}                               % Place acknowledgements
% Partially supported by the Roman Senate.  % here.
% \end{ack}

\appendix

\subsection{Proof of Lemma~\ref{lemma:bound_on_g}}    % Each appendix must have a short title.

For (1), by the projection's non-expansive property, we have
\begin{align*}
\bigg\|\Pi_{\mathcal{X}}\bigg[\sum_{j=1}^N [\mathbf{W}_r]_{ij}\mathbf{x}^j_t + \delta \mathbf{y}^i_t - \gamma_t g^t_{\mu^i}(\mathbf{x}^i_t)\bigg]-&\sum_{j=1}^N [\mathbf{W}_r]_{ij}\mathbf{x}^j_t\bigg\|\\
&\quad\leq \|\delta \mathbf{y}^i_t - \gamma_t g^t_{\mu^i}(\mathbf{x}^i_t)\|,
\end{align*}
which implies
\begin{align}
\|\theta^i_t\|&\leq \bigg\|\mathbf{x}^i_{t+1}-\sum_{j=1}^N [\mathbf{W}_r]_{ij}\mathbf{x}^j_t\bigg\|+\|\delta \mathbf{y}^i_t\|\nonumber\\
&\leq \|\delta \mathbf{y}^i_t - \gamma_t g^t_{\mu^i}(\mathbf{x}^i_t)\|+\delta \|\mathbf{y}^i_t\|\nonumber\\
&\leq\gamma_t\|g^t_{\mu^i}(\mathbf{x}^i_t)\|+2\delta \|\mathbf{y}^i_t\|.\label{eq:gik_ineq}
\end{align}
For the boundedness of $\|\mathbf{y}^i_t\|$, it follows from \eqref{eq:d-dgd} that for $t \geq 1$ 
\begin{align*}
\phi^i_t &= \sum_{j=1}^{2N}[\mathbf{W}^t]_{ij}\phi^j_0 + \sum_{r=1}^{t-1}\sum_{j=1}^{2N}[\mathbf{W}^{t-r}]_{ij}\theta^j_{r-1} + \theta^i_{t-1}.
\end{align*}
Thus, for $i = \{N+1,\ldots,2N\}$ and $t\geq1$ we have
\begin{align*}
&\|\phi^i_t\| \leq \sum_{j=1}^{2N}\big\|[\mathbf{W}^t]_{ij}\big\|\rho+ \sum_{r=1}^{t-1}\sum_{j=1}^{N}\big\|[\mathbf{W}^{t-r}]_{ij}\big\|\|\theta^j_{r-1}\|.
\end{align*}
Applying Lemma~2 yields
\begin{align*}
&\|\mathbf{y}^i_t\| \leq 2N\rho C\lambda^t+  C\sum_{r=1}^{t-1}\lambda^{t-r}\Theta_{r-1}.
\end{align*}
Thus, with the above result, it can be obtained from \eqref{eq:gik_ineq} that
\begin{align}
\Theta_t &\leq 4N^2\rho\delta C\lambda^t + \gamma_t\sum_{i=1}^{N}\|g^t_{\mu^i}(\mathbf{x}^i_t)\|+2N\delta C\sum_{r=1}^{t-1}\lambda^{t-r}\Theta_{r-1},\label{eq:gik_ineq_1}
\end{align}
Taking the total expectation and summing over $t = 1$ to $T$,
\begin{align*}
\sum_{t=1}^T\mathbf{E}[\Theta_t] &\leq \frac{4N^2\rho\delta C\lambda}{1-\lambda}+ N(p+4)\hat{D}\sum_{t=1}^T\gamma_t+\frac{2N\delta C\lambda}{1-\lambda}\sum_{t=1}^T\mathbf{E}[\Theta_t].
\end{align*}
The result (1) follows by moving the last term to the left hand side, adding the extra constant term for $t=0$ on both sides, and noting that $\delta < \frac{1-\lambda}{2\sqrt{3}NC\lambda}<\frac{1-\lambda}{2NC\lambda}$.

For (2), squaring both sides of \eqref{eq:gik_ineq_1}, we obtain
\begin{align*}
\Theta^2_t &\leq 48N^4\rho^2\delta^2 C^2\lambda^{2t} + 3N\gamma^2_t\sum_{i=1}^{N}\|g^t_{\mu^i}(\mathbf{x}^i_t)\|^2 +12N^2\delta^2 C^2\bigg(\sum_{r=1}^{t-1}\lambda^{t-r}\Theta_{r-1}\bigg)^2.
\end{align*}
Applying Cauchy-Schwarz inequality gives
\begin{align*}
\bigg(\sum_{r=1}^{t-1}\lambda^{t-r}\Theta_{r-1}\bigg)^2 &\leq \bigg(\sum_{r=1}^{t-1}\lambda^{t-r}\bigg)\bigg(\sum_{r=1}^{t-1}\lambda^{t-r}\Theta^2_{r-1}\bigg)\leq\frac{\lambda}{1-\lambda}\sum_{r=1}^{t-1}\lambda^{t-r}\Theta^2_{r-1}.
\end{align*}
Combining the above two inequalities and taking the total expectation, we have
\begin{align}
\mathbf{E}[\Theta^2_t] &\leq 48N^4\rho^2\delta^2 C^2\lambda^{2t} + 3N\gamma^2_t\sum_{i=1}^{N}\mathbf{E}[\|g^t_{\mu^i}(\mathbf{x}^i_t)\|^2] +\frac{12N^2\delta^2 C^2\lambda}{1-\lambda}\sum_{r=1}^{t-1}\lambda^{t-r}\mathbf{E}[\Theta^2_{r-1}].\label{eq:gik_square_ineq}
\end{align}
Invoking Lemma~\ref{lemma:property_f_mu}-3) and summing over time from $t = 1$ to $T$,
\begin{align*}
\sum_{t=1}^T\mathbf{E}[\Theta^2_t] &\leq 48N^4\rho^2\delta^2 C^2\sum_{t=1}^T\lambda^{2t}+ 3N^2(p+4)^2\hat{D}^2\sum_{t=1}^T\gamma^2_t +\frac{12N^2\delta^2C^2\lambda^2}{(1-\lambda)^2}\sum_{t=1}^T\mathbf{E}[\Theta^2_t],
\end{align*}
where we have applied 
\begin{align*}
\sum_{r=1}^{t-1}\lambda^{t-r}\mathbf{E}[\Theta^2_{r-1}]\leq\frac{\lambda}{1-\lambda}\sum_{t=1}^T\mathbf{E}[\Theta^2_t].
\end{align*}
Thus, the result (2) follows by moving the last term to the left hand side, adding the extra constant term for $t=0$ on both sides, and noting that $\delta < \frac{1-\lambda}{2\sqrt{3}NC\lambda}$.

For (3), taking square root on both sides of \eqref{eq:gik_square_ineq}, and summing over time from $t = 1$ to $T$
\begin{align*}
\sum_{t=1}^T\sqrt{\mathbf{E}[\Theta^2_t]} &\leq 4\sqrt{3}N^2\rho\delta C\sum_{t=1}^T\lambda^t + \sqrt{3}N(p+4)\hat{D}\sum_{t=1}^T\gamma_t+\frac{2\sqrt{3}N\delta C\lambda}{1-\lambda}\sum_{t=1}^T\sqrt{\mathbf{E}[\Theta^2_t]},
\end{align*}
where we have applied $\sqrt{\sum_i a^2_i} \leq \sum_i a_i$ ($a_i\geq 0$) and
\begin{align*}
\sum_{t=1}^T\sum_{r=1}^{t-1}\sqrt{\lambda}^{t-r}\sqrt{\mathbf{E}[\Theta^2_{r-1}]}\leq \frac{\sqrt{\lambda}}{\sqrt{1-\lambda}}\sum_{t=1}^T\sqrt{\mathbf{E}[\Theta^2_t]},
\end{align*}
The result (3) follows by moving the last term to the left hand side, adding the extra constant term for $t=0$ on both sides, and noting that $\delta < \frac{1-\lambda}{2\sqrt{3}NC\lambda}$.

\bibliographystyle{IEEEtran}
\bibliography{rgf_d_dpgd_time-var_reference}
                                 % bibliography (preferred). The
                                 % correct style is generated by
                                 % Elsevier at the time of printing.

\end{document}